\documentstyle[12pt]{article}
\newcommand{\sect}[1]{\section{#1}\setcounter{equation}{0}}

\font\mbn=msbm10 scaled \magstep1
\font\mbs=msbm7 scaled \magstep1
\font\mbss=msbm5 scaled \magstep1
\newfam\mbff
\textfont\mbff=\mbn
\scriptfont\mbff=\mbs
\scriptscriptfont\mbff=\mbss\def\mbf{\fam\mbff}

\def\Co{{\mbf C}}

\def\N{{\mbf N}}
\newtheorem{Th}{Theorem}[section]
\newtheorem{Lm}[Th]{Lemma}

\newtheorem{Proposition}[Th]{Proposition}

\author{Alexander Brudnyi
\thanks{Communicated by K. Saito. Received November 28, 2005.\newline
\newline
Department of Mathematics and Statistics,
University of Calgary, 2500 University Drive N.W., Calgary, Alberta
T2N 1N4, Canada. E-mail: albru@math.ucalgary.ca.
\newline
Research supported in part by NSERC.}}
\title{On Holomorphic $L_{2}$ functions on Coverings of Strongly
Pseudoconvex Manifolds}
\date{} 
\begin{document} 
\maketitle
\begin{abstract}
{In this paper we answer an important question posed in the paper [GHS] by
Gromov, Henkin and Shubin on existence of sufficiently many
holomorphic $L_{2}$ functions on arbitrary coverings of strongly pseudoconvex 
manifolds.} 

\end{abstract}
\sect{\hspace*{-1em}. Introduction.}
{\bf 1.1.} Let $M\subset\subset N$ be a domain with smooth boundary 
$bM$ in an $n$-dimensional complex manifold $N$, specifically,
\begin{equation}\label{m1}
M=\{z\in N\ :\ \rho(z)<0\}
\end{equation}
where $\rho$ is a real-valued function of class $C^{2}(\Omega)$ in a
neighbourhood $\Omega$ of the compact set $\overline{M}:=M\cup bM$
such that
\begin{equation}\label{m2}
d\rho(z)\neq 0\ \ \ {\rm for\ all}\ \ \ z\in bM\ .
\end{equation}
Let $z_{1},\dots, z_{n}$ be complex local coordinates in $N$ near $z\in bM$.
Then the tangent space $T_{z}N$ at $z$ is identified with $\Co^{n}$.
By $T_{z}^{c}(bM)\subset T_{z}N$ we denote the complex tangent space to
$bM$ at $z$, i.e.,
\begin{equation}\label{m3}
T_{z}^{c}(bM)=\{w=(w_{1},\dots,w_{n})\in T_{z}(N)\ :\ \sum_{j=1}^{n}
\frac{\partial\rho}{\partial z_{j}}(z)w_{j}=0\}\ .
\end{equation}
The {\em Levi form} of $\rho$ at $z\in bM$ is a hermitian form on
$T_{z}^{c}(bM)$ defined in the local coordinates by the formula
\begin{equation}\label{m4}
L_{z}(w,\overline{w})=\sum_{j,k=1}^{n}
\frac{\partial^{2}\rho}{\partial z_{j}\partial\overline{z}_{k}}(z)w_{j}
\overline{w}_{k}\ .
\end{equation}
The manifold $M$ is called {\em pseudoconvex} if $L_{z}(w,\overline{w})\geq 0$
for all $z\in bM$ and $w\in T_{z}^{c}(bM)$. It is called {\em strongly
pseudoconvex} if $L_{z}(w,\overline{w})>0$ for all $z\in bM$ and all
$w\neq 0$, $w\in T_{z}^{c}(bM)$.

Equivalently, strongly pseudoconvex manifolds can be described as the ones 
which locally, in a neighbourhood of any boundary point, can be presented as 
strictly convex domains in $\Co^{n}$. It is also known (see [C], [R]) that 
any strongly pseudoconvex manifold admits a proper holomorphic map with 
connected fibres onto a normal Stein space.

Without loss of generality we may and will assume that 
$\pi_{1}(M)=\pi_{1}(N)$ for $M$ as above.
Let $r: N'\to N$ be an unbranched covering of $N$. By
$M':=r^{-1}(M)$ we denote the corresponding covering of $M$. 
Also,  by $bM':=r^{-1}(bM)$ and $\overline{M'}:=M'\cup bM'$ we denote the 
boundary and the closure of $M'$ in $N'$. 

Let $dV_{M'}$ be the Riemannian volume form on $M'$ obtained by
a Riemannian metric pulled back from $N$.
By $H^{2}(M')$ we denote the
Hilbert space of holomorphic functions $g$ on $M'$ with norm
\begin{equation}\label{e1.5}
\left(\int_{z\in M'}|g(z)|^{2} dV_{M'}(z)\right)
^{1/2} .
\end{equation}

Let $X$ be a subspace of the space ${\cal O}(M')$ of all holomorphic
functions on $M'$.

A point $z\in bM'$ is called a {\em peak point} for $X$ if there 
exists a function $f\in X$ such that $f$ is unbounded on 
$M'$ but bounded outside $U\cap M'$ for any neighbourhood $U$ of $z$ in $N'$.

A point $z\in bM'$ is called a {\em local peak point} for 
$X$ if there exists a function $f\in X$ such that $f$ is unbounded in 
$U\cap M'$ for any neighbourhood $U$ of $z$ in $N'$ and there exists a 
neighbourhood $U$ of $z$ in $N'$ such that for any neighbourhood $V$ 
of $z$ in $N'$ the function $f$ is bounded on $U\setminus V$.

The Oka-Grauert theorem [G] implies that if $M$ is strongly 
pseudoconvex and $bM$ is not empty then every $z\in bM$ is a peak point for 
$H^{2}(M)$. In general it
is not known whether a similar statement is true for boundary points of 
an infinite covering $M'$ of $M$.

Assume that $M':=M_{G}$ is a regular covering of $M$ with a transformation
group $G$. In [GHS] the von Neumann $G$-dimension $dim_{G}$ was used to 
measure the space $H^{2}(M_{G})$. In particular, one of the
main results [GHS, Th.$\!$ 0.2] states:\\
\\
{\bf Theorem }\ {\em If $M$ is strongly pseudoconvex, then}
\begin{itemize}
\item[(a)]
$dim_{G}H^{2}(M_{G})=\infty$\ \ \ {\em and}
\item[(b)]
{\em each point in $bM_{G}$ is a local peak point for $H^{2}(M_{G})$.}
\end{itemize}

Also, in [GHS, p.$\!$ 3] the following important question was asked: ''{\em A natural
question arises: is the cocompact group action {\rm (on $\overline{M'}$)} really relevant for
the existence of many holomorphic $L^{2}$-functions {\rm (on $M'$)} or is it just an artifact of the chosen
methods which require a use of von Neumann algebras?} '' And further: ''{\em It is
not clear how to formulate conditions assuring that $dim\ \!L^{2}{\cal O}(M')=\infty$\footnote{$dim\ \!H^{2}(M')$ in our notation} without any group
action.}'' 

The main result of the present paper answers the above formulated questions. 
In particular, we show that the regularity of $M'$ is irrelevant for
the existence of many holomorphic $L^{2}$-functions on $M'$.  
Moreover, we also prove a substantial extension of the
above result of [GHS]. Our method of the proof is completely different
(and probably more easier) 
from that of used in [GHS] and is based on the $L_{2}$-cohomology 
techniques, as well as, on the geometric properties of $M$. \\
{\bf 1.2.} To formulate our result, let $C_{M}\subset M$ be the union of
all compact complex subvarieties of $M$ of 
complex dimension $\geq 1$. It is known that if $M$ is strongly
pseudoconvex, then $C_{M}$ is a compact complex subvariety of $M$.
Let $z_{i}$, $1\leq i\leq m$, be distinct points in $M\setminus C_{M}$.
By $l_{2}(z_{i}')$ we denote the Hilbert space of $l_{2}$ functions on the
fibre $z_{i}':=r^{-1}(z_{i})$. 
\begin{Th}\label{te1}
If $M$ is strongly pseudoconvex, then 
\begin{itemize}
\item[(a)]
For any $f_{i}\in l_{2}(z_{i}')$, $1\leq i\leq m$, there exists 
$F\in H^{2}(M')$
such that $F|_{z_{i}'}=f_{i}$, $1\leq i\leq m$;
\item[(b)]
Each point in $bM'$ is a peak point for $H^{2}(M')$. 
\end{itemize}
\end{Th}

Similar results are valid for certain $L_{p}$ spaces of 
holomorphic functions on $M'$. These and some other results
will be published elsewhere. It is worth noting that results
much stronger than Theorem \ref{te1} can be obtained if
$M$ is a strongly pseudoconvex Stein manifold, see [Br1], [Br2] for an 
exposition.
\sect{\hspace*{-1em}. Auxiliary Results.}
{\bf 2.1.}
Let $X$ be a complete K\"{a}hler manifold of dimension $n$ with a 
K\"{a}hler form $\omega$ and
$E$ be a hermitian holomorphic vector bundle on $X$ with curvature $\Theta$. 
Let $L_{2}^{p,q}(X,E)$ be the space of $L_{2}$ $E$-valued $(p,q)$-forms on $X$
with the $L_{2}$ norm, and let $W_{2}^{p,q}(X,E)$ be the subspace of forms
such that $\overline\partial\eta$ is $L_{2}$. (The forms $\eta$ may be taken 
to be either smooth or just measurable, in which case 
$\overline\partial\eta$ is understood in the distributional sense.) The
cohomology of the resulting $L_{2}$ Dolbeault complex 
$(W_{2}^{\cdot ,\cdot},\overline\partial)$ is the $L_{2}$ cohomology
$$
H_{(2)}^{p,q}(X,E)=Z_{2}^{p,q}(X,E)/B_{2}^{p,q}(X,E)\ ,
$$
where $Z_{2}^{p,q}(X,E)$ and $B_{2}^{p,q}(X,E)$ are the spaces of 
$\overline\partial$-closed and $\overline\partial$-exact forms in
$L_{2}^{p,q}(X,E)$, respectively.

If $\Theta\geq\epsilon\omega$
for some $\epsilon>0$ in the sense of Nakano, then the $L_{2}$ 
Kodaira-Nakano vanishing theorem, see [D], [O], states that
\begin{equation}\label{e1}
H_{(2)}^{n,r}(X,E)=0\ \ \ {\rm for}\ \ \ r>0\ .
\end{equation}
{\bf 2.2.}
Let $M\subset\subset N$ be a strongly pseudoconvex manifold. Without
loss of generality we will assume that $\pi_{1}(M)=\pi_{1}(N)$ and
$N$ is strongly pseudoconvex, as well. Then there exist a normal
Stein space $X_{N}$, a proper holomorphic surjective map $p:N\to X_{N}$ with
connected fibres and points $x_{1},\dots, x_{l}\in X_{N}$ such that
$$
p:N\setminus\bigcup_{1\leq i\leq l}p^{-1}(x_{i})\to 
X_{N}\setminus\bigcup_{1\leq i\leq l}\{x_{i}\}
$$ 
is biholomorphic, see [C], [R]. By definition, the domain 
$X_{M}:=p(M)\subset X_{N}$ is strongly pseudoconvex, and so it is
Stein. Without loss of generality we may assume that 
$x_{1},\dots, x_{l}\in X_{M}$. Thus 
$\cup_{1\leq i\leq l}\ p^{-1}(x_{i})=C_{M}$. 

Let $L\subset\subset N$ be a strongly pseudoconvex 
neighbourhood of $\overline{M}$. Then $X_{L}:=p(L)$ is a
strongly pseudoconvex neighbourhood of $\overline{X}_{M}$ in $X_{N}$.
We introduce a complete K\"{a}hler
metric on the complex manifold $L\setminus C_{M}$ as follows.

According to [N] there is a proper one-to-one holomorphic map
$i:X_{L}\hookrightarrow\Co^{2n+1}$, $n=dim_{\Co}X_{L}$, which is an embedding at regular points of $X_{L}$. Thus $i(X_{L})\subset\Co^{2n+1}$ is a 
closed complex subvariety. By $\omega_{e}$ we denote
the $(1,1)$-form on $L$ obtained as the pullback by 
$i\circ p$ of the Euclidean K\"{a}hler form on $\Co^{2n+1}$. Clearly,
$\omega_{e}$ is $d$-closed and positive outside $C_{M}$.

Similarly we can embed $X_{N}$ into $\Co^{2n+1}$ as a closed
complex subvariety. Let $j:X_{N}\hookrightarrow\Co^{2n+1}$ be an
embedding such that $j(X_{L})$ belongs to the open Euclidean ball $B$ of 
radius $1/4$ centered at $0\in\Co^{2n+1}$. 
Set $z_{i}:=j(x_{i})$, $1\leq i\leq l$. By
$\omega_{i}$ we denote the restriction to $L\setminus C_{M}$ of
the pullback with respect to $j\circ p$ of the form
$-\sqrt{-1}\cdot\partial\overline\partial\log(\log ||z-z_{i}||^{2})^{2}$
on $\Co^{2n+1}\setminus\{z_{i}\}$. (Here $||\cdot||$ stands for
the Euclidean norm on $\Co^{2n+1}$.) Since $j(X_{L})\subset B$, 
the form $\omega_{i}$ is K\"{a}hler.
Its positivity follows from the fact that the 
function $-\log(\log||z||^{2})^{2}$ is 
strictly plurisubharmonic for $||z||<1$.
Also, $\omega_{i}$ is extended to a smooth form on $L\setminus p^{-1}(x_{i})$.
Now, let us introduce a K\"{a}hler
form $\omega_{L}$ on $L\setminus C_{M}$ by the formula
\begin{equation}\label{e2}
\omega_{L}:=\omega_{e}+\sum_{1\leq i\leq l}\omega_{i}\ .
\end{equation}
\begin{Proposition}\label{p1}
The path metric $d$ on $L\setminus C_{M}$ induced by $\omega_{L}$ is complete.
\end{Proposition}
{\bf Proof.} 
Assume, on the contrary, that there is a sequence $\{w_{j}\}$
convergent either to $C_{M}$ or to the boundary $bL$ of $L$ such that the
sequence $\{d(o,w_{j})\}$ is bounded (for a fixed point 
$o\in L\setminus C_{M}$). Then, since $\omega_{L}\geq\omega_{e}$, 
the sequence $\{i(p(w_{j}))\}\subset\Co^{2n+1}$ is
bounded. This implies that $\{w_{j}\}$ converges to $C_{M}$.
But since $\omega_{L}\geq\sum\omega_{i}$, the latter is impossible, 
see, e.g., [GM] for similar arguments.\ \ \ \ \ $\Box$

In the same way one obtains complete K\"{a}hler metrics on unbranched 
coverings of $L\setminus C_{M}$ induced by
pullbacks to these coverings of the K\"{a}hler form $\omega_{L}$ 
on $L\setminus C_{M}$.\\
{\bf 2.3.} We retain the notation of the previous section. Also,
for an $n$-dimensional complex manifold $X$ by $T_{X}$ and $T^{*}_{X}$ we 
denote complex tangent and cotangent bundles on $X$ and 
by $K_{X}=\wedge^{n} T_{X}^{*}$ the canonical line bundle on $X$.


Let $r:N'\to N$ be an unbranched covering. Consider the corresponding
covering $(L\setminus C_{M})':=r^{-1}(L\setminus C_{M})$ of 
$L\setminus C_{M}$. 
We equip $(L\setminus C_{M})'$ with the complete
K\"{a}hler metric induced by the form
$\omega_{L}':=r^{*}\omega_{L}$. 

Next we consider the function $f:=\sum_{0\leq i\leq l}f_{i}$ 
on $(L\setminus C_{M})'$ such that
$f_{0}$ is the pullback by $i\circ p\circ r$ of 
the function $||z||^{2}$ on $\Co^{2n+1}$ and 
$f_{i}$ is the pullback by
$j\circ p\circ r$ of the function $-\log(\log ||z-z_{i}||^{2})^{2}$ on
$\Co^{2n+1}\setminus\{z_{i}\}$, $1\leq i\leq l$.
Clearly we have
\begin{equation}\label{e6}
\omega_{L}':=\sqrt{-1}\cdot\partial\overline\partial f\ .
\end{equation}

Let $E:=(L\setminus C_{M})'\times\Co$ be the trivial holomorphic line bundle
on $(L\setminus C_{M})'$. Let $g$ be the pullback to $(L\setminus C_{M})'$ of
a smooth plurisubharmonic function on $L$. 
We equip $E$ with the hermitian metric $e^{-f-g}$
(i.e., for $z\times v\in E$ the square of its norm in this
metric equals $e^{-f(z)-g(z)}|v|^{2}$ where $|v|$ is the modulus of 
$v\in\Co$). Then the curvature $\Theta_{E}$ of $E$ satisfies
\begin{equation}\label{e7}
\Theta_{E}:=-\sqrt{-1}\cdot\partial\overline\partial
\log(e^{-f-g})
=\omega_{L}'+\sqrt{-1}\cdot
\partial\overline\partial g\geq\omega_{L}'.
\end{equation}
Thus we can apply the $L_{2}$ Kodaira-Nakano vanishing theorem of section
2.1 to get
\begin{equation}\label{e9}
H_{(2)}^{n,r}((L\setminus C_{M})',E)
=0\ \ \ {\rm for}\ \ \ r>0\ .
\end{equation}

Let $K_{(L\setminus C_{M})'}$ be the canonical holomorphic line bundle on
$(L\setminus C_{M})'$ equipped with the hermitian metric induced by
$\omega_{L}'$. Consider the hermitian line bundle 
$V_{g}:=E\otimes 
K_{(L\setminus C_{M})'}$ equipped with the tensor product of the 
corresponding hermitian metrics. Then from (\ref{e9}) we have
\begin{equation}\label{e10}
H_{(2)}^{0,r}((L\setminus C_{M})', V_{g})\cong 
H_{(2)}^{n,r}((L\setminus C_{M})',E)
=0\ \ \ {\rm for}\ \ \ r>0\ .
\end{equation}
{\bf 2.4.} Let $U\subset L$ be a relatively compact neighbourhood of
$C_{M}$. Consider a finite open cover $(U_{i})_{1\leq i\leq k}$ of
$\overline{L\setminus U}$ by simply connected coordinate charts
$U_{i}\subset\subset N\setminus C_{M}$. 
We introduce complex coordinates on  
$U_{i}':=r^{-1}(U_{i})\subset N'$ by the pullback of the coordinates on 
$U_{i}$. In these coordinates $U_{i}'$ is naturally identified with 
$U_{i}\times S$ where $S$ is the fibre of $r:N'\to N$. 

Let $\eta$ be a smooth $(p,q)$-form on $(L\setminus C_{M})'$ equals 0 on 
$r^{-1}(U)$. Then in the above holomorphic coordinates $(z,s)$,
$z=(z_{1},\dots, z_{n})\in U_{i}\cap L$, $s\in S$, on $U_{i}'\cap L'$,
$L':=r^{-1}(L)$, the form $\eta$ is presented as
\begin{equation}\label{e11}
\eta(z,s)=\sum_{i_{1},\dots, i_{p},j_{1},\dots,j_{q}}
\eta_{i;i_{1},\dots, i_{p},j_{1},\dots,j_{q}}(z,s)\ \!dz_{i_{1}}
\wedge\cdots\wedge
dz_{i_{p}}\wedge d\overline{z}_{j_{1}}\wedge\cdots\wedge d\overline{z}_{j_{q}}
\end{equation}
where $\eta_{i;i_{1},\dots, i_{p},j_{1},\dots,j_{q}}$ are smooth functions on
$(U_{i}\cap L)\times S$. 

We say that $\eta$ belongs to the space
${\cal E}_{U;2}^{p,q}((L\setminus C_{M})')$ if in (\ref{e11}) we have
\begin{equation}\label{e12}
\sup_{z\in U_{i}\cap L,i,i_{1},\dots, i_{p},j_{1},\dots,j_{q}}
\left\{\sum_{s\in S}
|\eta_{i;i_{1},\dots, i_{p},j_{1},\dots,j_{q}}(z,s)|^{2}
\right\}<\infty\ .
\end{equation}

Let $e$ be a holomorphic section of $K|_{L\setminus C_{M}}$. Then
$\eta\cdot r^{*}e$ can be viewed as a $(p,q)$-form with values in $V_{g}$. 
(Here $r^{*}e$ is the pullback of $e$ to $(L\setminus C_{M})'$, i.e., 
$r^{*}e\in {\cal O}((L\setminus C_{M})',K|_{(L\setminus C_{M})'})$.)
\begin{Proposition}\label{pr2}
For every $\eta\in {\cal E}_{U;2}^{p,q}((L\setminus C_{M})')$ and
$e\in {\cal O}(L\setminus C_{M},K|_{L\setminus C_{M}})$
there is a plurisubharmonic function $g$ in the definition of the 
metric on $V_{g}$ such that $\eta\cdot r^{*}e\in 
L_{2}^{p,q}((L\setminus C_{M})', V_{g})$.
\end{Proposition}
{\bf Proof.} 
In this proof by $||\cdot||$ we denote the hermitian metric on 
the space of $V_{g}$-valued $(p,q)$-forms induced by the hermitian metrics on
$V_{g}$ and $T_{(L\setminus C_{M})'}$. Set
$$
h_{i;i_{1},\dots, i_{p},j_{1},\dots,j_{q}}(z):=||r^{*}e(z,s)\cdot dz_{i_{1}}
\wedge\cdots\wedge
dz_{i_{p}}\wedge d\overline{z}_{j_{1}}\wedge\cdots\wedge 
d\overline{z}_{j_{q}}||^{2}
$$
Then $h_{i;i_{1},\dots, i_{p},j_{1},\dots,j_{q}}$ is a nonnegative 
continuous function on $U_{i}\cap L$. Let $\hat g$ be such that
$r^{*}\hat g=g$. By the definition of metrics
on $V_{g}$ and $T_{(L\setminus C_{M})'}$ 
\begin{equation}\label{e314}
h_{i;i_{1},\dots, i_{p},j_{1},\dots,j_{q}}(z):=
\hat h_{i;i_{1},\dots, i_{p},j_{1},\dots,j_{q}}(z)\cdot 
e^{-\hat g(z)},
\end{equation}
where $\hat h_{i;i_{1},\dots, i_{p},j_{1},\dots,j_{q}}$ is a 
nonnegative continuous function on $U_{i}\cap L$ independent of 
$\hat g$. 

Now for some $A\in\N$ we have
\begin{equation}\label{e13}
||\eta(z,s)\cdot r^{*}e(z,s)||^{2}\leq
A\times\sum_{i_{1},\dots, i_{p},j_{1},\dots,j_{q}}
|\eta_{i;i_{1},\dots, i_{p},j_{1},\dots,j_{q}}(z,s)|^{2}\cdot
h_{i;i_{1},\dots, i_{p},j_{1},\dots,j_{q}}(z)\ .
\end{equation}
According to the definition of $L_{2}^{p,q}((L\setminus C_{M})', V_{g})$ 
we have to show that
$$
|\eta\cdot r^{*}e|^{2}:=\int_{(L\setminus C_{M})'}||\eta\cdot r^{*}e||^{2}\cdot 
(\omega_{L}')^{n}<\infty\ .
$$
Since $\omega_{L}'=r^{*}\omega_{L}$, from (\ref{e314}) and (\ref{e13}) we get
$$
\begin{array}{l}
|\eta\cdot r^{*}e|^{2}\leq\\
\\
\displaystyle
A\times\sum_{i=1}^{k}\int_{U_{i}\cap L}\left(
\sum_{i_{1},\dots, i_{p},j_{1},\dots,j_{q},
s\in S}|\eta_{i;i_{1},\dots, i_{p},j_{1},\dots,j_{q}}(\cdot,s)|^{2}\right)
\hat h_{i;i_{1},\dots, i_{p},j_{1},\dots,j_{q}}e^{-\hat g}
\omega_{L}^{n}.
\end{array}
$$
Also, by the hypothesis of the proposition, see (\ref{e12}), 
$$
\sup_{z\in U_{i}\cap L}\left\{\sum_{i_{1},\dots, i_{p},j_{1},\dots,j_{q},
s\in S}|\eta_{i;i_{1},\dots, i_{p},j_{1},\dots,j_{q}}(z,s)|^{2}
\right\}
<\infty\ \ \
{\rm for}\ \ \ 1\leq i\leq k\ .
$$
Thus in order to prove the proposition it suffices to check that there is 
$\hat g$ in the definition of the metric on $V_{g}$ such that for every $i$
\begin{equation}\label{e14}
\int_{U_{i}\cap L}
\hat h_{i;i_{1},\dots, i_{p},j_{1},\dots,j_{q}}e^{-\hat g}
\omega_{L}^{n}<\infty\ .
\end{equation}
The required result now follows from 
\begin{Lm}\label{le1}
Let $h$ be a nonnegative piecewise continuous function on $L$ equals 0 in 
some neighbourhood of $C_{M}$ and bounded on every compact subset of
$L\setminus C_{M}$. Then there exists a smooth
plurisubharmonic function $\hat g$ on $L$ such that 
$$
\int_{L}h\cdot e^{-\hat g}\ \!\omega_{L}^{n}<\infty\ .
$$
\end{Lm}
{\bf Proof.} Without loss of generality we identify $L\setminus C_{M}$
with $X_{L}\setminus\cup_{1\leq j\leq l}\ \{x_{j}\}$. Also, we identify
$X_{L}$ with a closed subvariety of $\Co^{2n+1}$ as in section 2.2. 
Let $U$ be a neighbourhood of $\cup_{1\leq j\leq l}\ \{x_{j}\}$ such that
$h|_{U}\equiv 0$. By $\Delta_{r}\subset\Co^{2n+1}$ we denote the open
polydisk 
of radius $r$ centered at $0\in\Co^{2n+1}$. Assume without loss of generality
that $0\in X_{L}\setminus U$. Consider the monotonically 
increasing function
\begin{equation}\label{e15}
v(r):=\int_{\Delta_{r}\cap (X_{L}\setminus U)}h\cdot\omega_{L}^{n}\ ,\ \ \
r\geq 0\ .
\end{equation}
By $v_{1}$ we denote a smooth monotonically increasing function satisfying
$v\leq v_{1}$ (such $v_{1}$ can be easily constructed
by $v$). Let us determine 
$$
v_{2}(r):=\int_{0}^{r+1}2v_{1}(2t)\ \!dt\ ,\ \ \ r\geq 0\ .
$$
By the definition $v_{2}$ is smooth, convex and monotonically
increasing. Moreover, 
$$
v_{2}(r)\geq\int_{\frac{r+1}{2}}^{r+1}2v_{1}(2t)\ \!dt\geq (r+1)v(r+1)\ .
$$
Next we define a smooth 
plurisubharmonic function $v_{3}$ on $\Co^{2n+1}$ by the formula
$$
v_{3}(z_{1},\dots, z_{2n+1}):=\sum_{j=1}^{2n+1}v_{2}(|z_{j}|)\ .
$$
Then the pullback of $v_{3}$ to $L$ is a smooth 
plurisubharmonic function on $L$. This is the required function 
$\hat g$. Indeed, under
the identification described at the beginning of the proof we have
$$
\begin{array}{c}
\displaystyle
\int_{L}h\cdot e^{-\hat g}\ \!\omega_{L}^{n}=\sum_{k=1}^{\infty}
\int_{(\Delta_{k}\setminus\Delta_{k-1})\cap (X_{L}\setminus U)}
h\cdot  e^{-\hat g}\ \!\omega_{L}^{n}\leq\\
\\
\displaystyle
\sum_{k=1}^{\infty}v(k)
e^{-v_{2}(k-1)}\leq\sum_{k=1}^{\infty}v(k)e^{-kv(k)}<\infty\ .\ \ \ \ \
\Box
\end{array}
$$

To complete the proof of the proposition it remains to put in the
above lemma
$$
h:=\sum_{i,i_{1},\dots, i_{p},j_{1},\dots,j_{q}}
\rho_{i}\cdot\hat h_{i;i_{1},\dots, i_{p},j_{1},\dots,j_{q}}
$$
where $\rho_{i}$ is the characteristic function of $U_{i}\cap L$.
\ \ \ \ \  $\Box$\\
{\bf 2.5.}  
Let $O\subset\subset L$ be a neighbourhood of $C_{M}$. We set
$O':=r^{-1}(O)$, $C_{M}':=r^{-1}(C_{M})$. Assume that the manifold $N$, see
section 1.1, is
equipped with a hermitian metric $\rho$. We equip the bundle $K_{L}$ with
the hermitian metric induced by $\rho$. Also, we equip 
$K_{L'}:=r^{*}K_{L}$ with the hermitian metric $\rho':=r^{*}\rho$.
\begin{Proposition}\label{pr3}
Any $h\in L_{2}((L\setminus C_{M})',V_{g})$ holomorphic
on $O'\setminus C_{M}'$ admits an extension to a section
$\hat h$ of $K_{L'}$ such that $\hat h|_{M'}\in L_{2}(M',K_{L'})$.
\end{Proposition}
{\bf Proof.} Consider a coordinate
neighbourhood $U\subset\subset O$ of a point $q\in C_{M}$ with coordinates
$z=(z_{1},\dots, z_{n})$. Taking the pullback of these
coordinates to $r^{-1}(U)$ we identify $r^{-1}(U)$ with 
$U\times S$ where $S$ is the fibre of $r$. Then 
$$
h(z,s)=h_{U}(z,s)dz_{1}\wedge\cdots\wedge dz_{n},\ \ \ z\in U\setminus C_{M},\ 
s\in S.
$$
By the definition of the metric $||\cdot||$ on $V_{g}$ we have
$$
\begin{array}{c}
\displaystyle
||h(z,s)||^{2}\omega_{L}^{n}(z)=
|h_{U}(z,s)|^{2}||dz_{1}\wedge\cdots\wedge dz_{n}||^{2}\omega_{L}^{n}(z)=\\
\\
\displaystyle
|h_{U}(z,s)|^{2}e^{-\hat g(z)}
(\sqrt{-1})^{n}\wedge_{i=1}^{n}dz_{i}\wedge d\overline{z}_{i}.
\end{array}
$$
Now, the hypotheses of the
proposition imply that
\begin{equation}\label{e16}
\int_{z\in U\setminus C_{M}}
\left(\sum_{s\in S}|h_{U}(z,s)|^{2}\right)e^{-\hat g(z)}
(\sqrt{-1})^{n}\wedge_{i=1}^{n}dz_{i}\wedge d\overline{z}_{i}<\infty\ .
\end{equation}
Let $\omega^{n}$ be the 
volume form induced by the hermitian metric $\rho$ on $N$ with the associated 
$(1,1)$-form $\omega$. Since by our construction $\hat g$ is smooth on $L$,
we have on $U$:
$$
e^{-\hat g(z)}(\sqrt{-1})^{n}\wedge_{i=1}^{n}dz_{i}\wedge d\overline{z}_{i}\geq
c\ \!\omega^{n}\ 
$$
for some $c>0$. From here and (\ref{e16}) we get
\begin{equation}\label{e17}
\int_{z\in U\setminus C_{M}}
\left(\sum_{s\in S}|h_{U}(z,s)|^{2}\right)\ \!
\omega^{n}(z)<\infty\ .
\end{equation}
In particular, this implies that 
every $h_{U}(\cdot, s)$, $s\in S$, belongs to the $L_{2}$ space
on $U\setminus C_{M}$ defined by integration with respect to the volume form
$(\sqrt{-1})^{n}\wedge_{i=1}^{n}dz_{i}\wedge d\overline{z}_{i}$. Also,
every $h_{U}(\cdot, s)$ is holomorphic on $U\setminus C_{M}$. Using these
facts and the Cauchy integral formulas for coefficients of the Laurent
expansion of $h_{U}(\cdot, s)$, one obtains easily that
every $h_{U}(\cdot, s)$ can be extended holomorphically to $U$. In turn, this
gives an extension $\hat h$ of $h$ to $r^{-1}(U)$. Now from (\ref{e17})
we obtain that $\hat h\in L_{2}(U,K_{L'})$.

Next assume that $\widetilde U\subset\subset L\setminus C_{M}$ is a  
simply connected coordinate neighbourhood of a point 
$q\in \overline{M}\setminus C_{M}$. Identifying
$r^{-1}(\widetilde U)$ with $\widetilde U\times S$ we have anew inequality 
of type (\ref{e16}) for $h|_{r^{-1}(\widetilde U)}$. 
Since  $\hat g$ is smooth on $L\setminus C_{M}$, repeating literally the
previous arguments we get that $h\in L_{2}(\widetilde U,K_{L'})$.
Taking a finite open cover of $\overline{M}$ by the above neighbourhoods $U$
and $\widetilde U$ and considering the extension of $h$ to
$M'$ defined by the above extended forms $\hat h$ on
$r^{-1}(U)$ we get the required result.\ \ \ \ \ $\Box$\\
\sect{\hspace*{-1em}. Proofs.}
{\bf 3.1. Proof of Theorem \ref{te1}$\ \!$(a).} 
First, we prove Theorem \ref{te1} (a) for $m=1$:
\begin{Th}\label{te3}
Let $z\in M\setminus C_{M}$ and $z':=r^{-1}(z)\in M'$. Then for any
$f\in l_{2}(z')$ there exists $F\in H^{2}(M')$ such that 
$F|_{z'}=f$.
\end{Th}
{\bf Proof.} In the proof we retain the notation of section 2.

Let $p:N\to X_{N}$ be the proper holomorphic map onto the 
normal Stein space $X_{N}$ from section 2.2 such that 
$p:N\setminus C_{M}\to X_{N}\setminus\cup_{1\leq i\leq l}\ \{x_{i}\}$ is
biholomorphic.
Since $X_{N}$ is Stein, there is a holomorphic function $h$ on 
$X_{N}$ whose set of zeros ${\cal Z}_{h}$
contains $p(z)$ and does not intersect $\cup_{1\leq j\leq l}\ \{x_{j}\}$. 
Let $O\subset\subset X_{N}$ be a Stein neighbourhood 
of the compact set ${\cal Z}_{h}\cap\overline{X}_{L}$, $X_{L}:=p(L)$, 
such that
$\overline{O}\cap\cup_{1\leq j\leq l}\ \{x_{j}\}=\emptyset$ and
$\overline{O}$ is holomorphically convex in $X_{N}$. 
We set
$O':=(p\circ r)^{-1}(O)\subset N'$. Then according to 
[Br1, Th. $\!$1.10] there is a holomorphic function $h_{1}$ on $O'$ satisfying
\begin{equation}\label{e31}
h_{1}|_{z'}=f\ \ \ {\rm and}\ \ \ 
\sup_{y\in p^{-1}(O)}\left\{\sum_{x\in r^{-1}(y)}|h_{1}(x)|^{2}
\right\}<\infty\ .
\end{equation}

Let $\rho$ be a $C^{\infty}$ function on 
$X_{N}\setminus\cup_{1\leq i\leq l}\ \{x_{i}\}$ equals $1$ in some
neighbourhood of ${\cal Z}_{h}\cap O$ in $O$ and $0$ outside $O$. 
By $\rho':=(p\circ r)^{*}\rho$ we denote its pullback to $N'$. Then 
$h_{2}:=(\rho'\cdot h_{1})|_{(L\setminus C_{M})'}$ is a $C^{\infty}$ 
function on $(L\setminus C_{M})'$. Let $h':=(p\circ r)^{*}h|_{L'}$ be the 
restriction to $L'$ of the pullback of $h$.  
Consider the $C^{\infty}$ $(0,1)$-from $\eta:=\overline\partial{h_{2}}/h'$.
It follows easily from (\ref{e31}) and (\ref{e12}) that 
$\eta\in {\cal E}_{U;2}^{0,1}((L\setminus C_{M})')$ for some
$U\subset L\setminus p^{-1}(O)$. 

Next, since $X_{N}$ can be embedded to $\Co^{2n+1}$, and since $z$ is a 
regular point of $X_{N}$, there is a section $e\in {\cal O}(N,K_{N})$ such that
$e(z)\neq 0$. Let $V\subset L$ be the set of zeros of $e$. Then $V$ is 
contained in the preimage $p^{-1}(V')$ of a complex analytic subspace
$V'\subset X_{N}$. Since $X_{N}$ is Stein, the latter impies that
there is a bounded holomorphic function $f_{e}$ on $L$ such that
$f_{e}\cdot e^{-1}$ is a bounded section from ${\cal O}(L,K_{L}^{*})$ and
$f_{e}(z)=1$ (recall that $L\subset\subset N$).  

Further, by Proposition \ref{pr2}
there is a function $g$ in the definition of the metric on the vector
bundle $V_{g}$ from section 2.3 such that 
$\widetilde\eta:=\eta\cdot r^{*}e\in L_{2}^{0,1}((L\setminus C_{M})',V_{g})$. 
Since 
$\overline\partial\widetilde\eta=0$, by (\ref{e10}) there is a
section $h_{3}\in L_{2}((L\setminus C_{M})',V_{g})$ such that
$\overline\partial h_{3}=\widetilde\eta$. Choose $g$ in the definition of the
metric on $V_{g}$ so that also
$h_{2}\cdot r^{*}e\in L_{2}((L\setminus C_{M})',V_{g})$. Thus 
$h_{4}:=h_{2}\cdot r^{*}e-h'h_{3}\in L_{2}((L\setminus C_{M})',V_{g})$, 
is holomorphic on $(L\setminus C_{M})'$. By Proposition \ref{pr3}
$h_{4}$ is extended to a holomorphic section $h_{4}'$ of $K_{L'}$ such
that $h_{4}'|_{M'}\in L_{2}(M',K_{L'})$. Moreover, by our construction
$h_{4}'|_{z'}=f\cdot r^{*}e|_{z'}$.
Finally, we set 
$$
F:=h_{4}'|_{M'}\cdot (r^{*}f_{e}\cdot r^{*}e^{-1})|_{M'}.
$$
Then $F\in H^{2}(M')$ and $F|_{z'}=f$. 
\ \ \ \ \ $\Box$

Now, let us prove Theorem \ref{te1} (a).

By Theorem \ref{te3} there are functions $F_{i}\in H^{2}(M')$
such that $F_{i}|_{z_{i}'}=f_{i}$, $1\leq i\leq m$. 
Let $t_{i}$ be a holomorphic
function on $N$ such that $t_{i}(z_{i})=1$ and $t_{i}(z_{j})=0$ for
$j\neq i$. Then the function $F:=\sum_{1\leq i\leq m} r^{*}t_{i}\cdot F_{i}$
satisfies the required condition.\ \ \ \ \ $\Box$\\
{\bf 3.2. Proof of Theorem \ref{te1}$\ \!$(b).} 
Let $x\in bM$ be a boundary point of
a strongly pseudoconvex manifold $M\subset\subset N$. We will prove that
each point $y\in r^{-1}(x)$ is a peak point for $H^{2}(M')$.

Choose a coordinate
neighbourhood $U_{1}\subset\subset N$ of $x$ with complex coordinates 
$w=(w_{1},\dots, w_{n})$ centered at $x$ such that the complex tangent space 
$T_{x}^{c}(bM)$ is given by the equation $w_{1}=0$. Consider the 
Taylor expansion of the defining function $\rho$ for $M$, see (\ref{m1}), 
(\ref{m2}), at $x$:
\begin{equation}\label{e32}
\rho(w)=\rho(x)+ 2 {\rm Re} f(x,w)+L_{x}(w-x,\overline{w}-\overline{x})+
O(||w-x||^{3})\ ,
\end{equation}
where $L_{x}$ is the Levi form at $x$ and $f(x,w)$ is a complex quadratic
polynomial with respect to $z$:
$$
f(x,w)=\frac{\partial\rho}{\partial w_{1}}(x)w_{1}+
\frac{1}{2}\sum_{1\leq\mu,\nu\leq n
}\frac{\partial^{2}\rho}{\partial w_{\mu}\partial w_{\nu}}(x)w_{\mu}w_{\nu}\ .
$$
(Here $\frac{\partial\rho}{\partial w_{1}}(x)\neq 0$ and 
$\frac{\partial\rho}{\partial w_{i}}(x)=0$ for $2\leq i\leq n$
by the choice of the coordinates on $U_{1}$). Next, in a sufficiently small
simply connected
neigbourhood $U\subset\subset U_{1}$ of $x$ we introduce new 
coordinates $z=(z_{1},\dots, z_{n})$ by the formulas
$$
z_{1}=f(x,w)\ \ \ {\rm and}\ \ \ z_{i}=w_{i}\ \ \ {\rm for}\ \ \
2\leq i\leq n\ .
$$
Since $bM$ is strongly pseudoconvex at $x$, diminishing if necessary $U$ we
have that in $U$ the intersection of the hyperplane $z_{1}=0$ with
$\overline{M}$ consists of one point $x$ and 
${\rm Re}\ \! z_{1}<0$ in $U$. Let $H:=\{z\in U\ :\ {\rm Re}\ \! z_{1}<0\}$.
Then we can choose a branch of $\log z_{1}$ so that 
$h_{x}(z)=\log z_{1}$ is a holomorphic function in $H$. Clearly,
$h_{x}\in L_{2}(H)$ for $L_{2}$ defined by a Riemannian
volume form pulled back from $N$, and $h_{x}$ has a peak point at $x$.

Let $H':=r^{-1}(H)\subset\subset N'$. Take a point 
$y\in r^{-1}(x)\subset bM'$. Let $H_{y}$ be the connected component of
$H'$ containing $y$. We set $h_{y}:=r^{*}h_{x}|_{H_{y}}$. Next, consider 
a smooth function $\rho$ on $N$ equals $1$ in a 
neigbourhood $O\subset\subset U$ of $x$ with support $S\subset\subset U$.
Let $O_{y}\subset\subset U_{y}$ be connected components of 
$O'=r^{-1}(O)$ and $U'=r^{-1}(U)$ containing
$y$. By $\rho_{y}$ we denote the pullback of $\rho|_{U}$ to $U_{y}$ and
by $S'\subset\subset U_{y}'$ support of $\rho_{y}$.
Then $\overline\partial (h_{y}\rho_{y})$ can be extended by $0$ to 
a $C^{\infty}$ $(0,1)$-form $\eta$ on 
$(N'\setminus S')\cup H_{y}'\cup O_{y}'$. Note that
$(N\setminus S)\cup H\cup O$ is a neighbourhood of $\overline M$. Hence, 
since $M\subset\subset N$ is strongly pseudoconvex, there is a strongly 
pseudoconvex manifold $L\subset\subset (N\setminus S)\cup H\cup O$ such that 
$\overline{M}\subset L$. Observe also that 
$r^{-1}((N\setminus S)\cup H\cup O)\subset 
(N'\setminus S')\cup H_{y}'\cup O_{y}'$. Thus the form $\eta$ is well-defined
on $L':=r^{-1}(L)$. Moreover, by our definitions
$\eta\in {\cal E}_{W;2}^{0,1}((L\setminus C_{M})')$ for some neighbourhood
$W$ of $C_{M}$, see section 2.4. 

As in the proof of part (a) we will choose a section $e\in {\cal O}(N,K_{N})$
such that $e(x)\neq 0$ and a function $f_{e}\in {\cal O}(L)$ such that
$f_{e}(x)=1$ and $f_{e}\cdot e^{-1}$ is a bounded section from 
${\cal O}(L,K_{L}^{*})$.

Now, by Proposition \ref{pr2}
there is a function $g$ in the definition of the metric on the vector
bundle $V_{g}$ from section 2.3 such that 
$\widetilde\eta:=\eta\cdot r^{*}e\in 
L_{2}^{0,1}((L\setminus C_{M})',V_{g})$.
Since $\overline\partial{\widetilde\eta}=0$, by (\ref{e10}) there is a section
$h_{1}\in L_{2}((L\setminus C_{M})',V_{g})$ such that 
$\overline\partial h_{1}=\widetilde\eta$. According to Proposition \ref{pr3}
$h_{1}|_{(M\setminus C_{M})'}$ is extended to a section $h_{2}$ of $K_{L'}$
such that $h_{2}\in L_{2}(M',K_{L'})$.
From here, using the fact that $(h_{y}\rho_{y})|_{M'\cap H_{y}}$ is
extended by $0$ to a smooth $L_{2}$ function on $M'$, we get
$h':=(h_{y}\rho_{y})\cdot r^{*}e-h_{2}\in H^{2}(M',K_{L'})$. Finally, we set
$$
h:=h'\cdot (r^{*}f_{e}\cdot r^{*}e^{-1})|_{M'}.
$$
Clearly, $h\in H^{2}(M')$ and it has a local peak 
point at $y$. Let us show that in fact $h$ has a peak 
point at $y$. This will complete the proof of the theorem.

Let $L_{1}\subset\subset L$ be a neighbourhood of $\overline M$. We set
$L_{1}':=r^{-1}(L_{1})$. Then by Proposition \ref{pr3} 
$h_{1}|_{r^{-1}(L_{1}\setminus C_{M})}$ is extended to a section of
$K_{L_{1}'}$ (denoted as before by $h_{2}$). By the definition, 
$h_{2}':=h_{2}\cdot (r^{*}f_{e}\cdot r^{*}e^{-1})|_{L_{1}'}$ is a
smooth holomorphic on $L_{1}'\setminus S'$ function. From the facts that
$L_{1}\subset\subset L$,
the $L_{2}$ norm $||\cdot||$ on $L_{1}'$ is defined 
by a Riemannian volume form pulled back from $N$ using the mean-value
property for the plurisubharmonic function $|h_{2}'|^{2}$ on 
$L_{1}'\setminus S'$ we get for some $c>0$,
$$
\sup_{z\in M}\left\{\sum_{w\in r^{-1}(z), w\notin U_{y}}
|h_{2}'(w)|^{2}\right\}\leq c\cdot ||h_{2}'||^{2}<\infty\ .
$$
In addition, $h_{2}'$ is continuous on $\overline{U}_{y}\cap\overline{M'}$.
These imply easily that $h$ is bounded outside $U\cap M'$ for any 
neighbourhood $U$ of $y$ in $N'$.\ \ \ \ \ $\Box$

\end{document}